\documentclass[11pt]{amsart}
\usepackage{amsthm,amsmath,graphicx,calc}
\usepackage[sort&compress,numbers]{natbib}
\usepackage[lmargin=35mm,rmargin=35mm,tmargin=30mm,bmargin=30mm]{geometry}
\renewcommand{\baselinestretch}{1.2}
\allowdisplaybreaks

\theoremstyle{plain}
\newtheorem{theorem}{Theorem}[section]
\newtheorem{lemma}[theorem]{Lemma}
\newtheorem{corollary}[theorem]{Corollary}
\newtheorem{proposition}[theorem]{Proposition}
\newtheorem{conjecture}[theorem]{Conjecture}

\theoremstyle{definition}

\newtheorem{open}[theorem]{Open Problem}

\newcommand{\thm}[2]{\begin{theorem}\label{thm:#1}#2\end{theorem}}
\newcommand{\con}[2]{\begin{conjecture}\label{con:#1}#2\end{conjecture}}
\newcommand{\cor}[2]{\begin{corollary}\label{cor:#1}#2\end{corollary}}
\newcommand{\lem}[3][]{\begin{lemma}[#1]\label{lem:#2}#3\end{lemma}}
\newcommand{\opn}[2]{\begin{open}\label{open:#1}#2\end{open}}
\newcommand{\prop}[2]{\begin{proposition}\label{prop:#1}#2\end{proposition}}
\newcommand{\prf}[2][Proof]{\begin{proof}[#1]#2\end{proof}}
\newcommand{\sect}[2]{\section{#2}\label{sec:#1}}

\newcommand{\corref}[1]{Corollary~\ref{cor:#1}}
\newcommand{\conref}[1]{Conjecture~\ref{con:#1}}
\newcommand{\propref}[1]{Proposition~\ref{prop:#1}}
\newcommand{\openref}[1]{Open Problem~\ref{open:#1}}
\newcommand{\thmref}[1]{Theorem~\ref{thm:#1}}
\newcommand{\lemref}[1]{Lemma~\ref{lem:#1}}

\newcommand{\threelemref}[3]{Lemmas~\ref{lem:#1}, \ref{lem:#2} and \ref{lem:#3}}
\newcommand{\secref}[1]{Section~\ref{sec:#1}}

\newcommand{\half}{\ensuremath{\protect\tfrac{1}{2}}}
\newcommand{\eighth}{\ensuremath{\protect\tfrac{1}{8}}}
\newcommand{\ninth}{\ensuremath{\protect\tfrac{1}{9}}}
\newcommand{\bracket}[1]{\ensuremath{\protect\left(#1\right)}}

\newcommand{\ceil}[1]{\ensuremath{\protect\lceil#1\rceil}}
\newcommand{\Oh}[1]{\ensuremath{\protect\mathcal{O}(#1)}}
\newcommand{\e}{\ensuremath{\mathbf{e}}}
\newcommand{\A}{\ensuremath{\mathcal{A}}}
\newcommand{\D}{\ensuremath{\mathcal{D}}}
\newcommand{\Z}{\ensuremath{\mathbb{Z}}}
\newcommand{\Prob}{\textbf{\textup{P}}}

\begin{document}

\title{Notes on Nonrepetitive Graph Colouring}

\thanks{The research of J.~Bar{\'a}t was supported by a Marie Curie
Fellowship of the European Community under contract number HPMF-CT-2002-01868
and by the OTKA Grant T.49398. The research of D.~Wood is supported by a Marie Curie Fellowship of the European Community under contract MEIF-CT-2006-023865, and by the projects MEC MTM2006-01267 and DURSI 2005SGR00692.}


\author{J{\'a}nos~Bar{\'a}t}
\address{\newline Department of Mathematics\newline University of Szeged\newline Szeged, Hungary}
\email{barat@math.u-szeged.hu}

\author{David~R.~Wood} 
\address{\newline Departament de Matem\'atica Aplicada II\newline 
Universitat Polit\`ecnica de Catalunya\newline
Barcelona, Spain}
\email{david.wood@upc.es}

\begin{abstract}
A vertex colouring of a graph is \emph{nonrepetitive on paths} if there is no path $v_1,v_2,\dots,v_{2t}$ such that $v_i$ and $v_{t+i}$ receive the same colour for all $i=1,2,\dots,t$. We determine the maximum density of a graph that admits a $k$-colouring that is nonrepetitive on paths. We prove that every graph has a subdivision that admits a $4$-colouring that is nonrepetitive on paths. The best previous bound was $5$. We also study colourings that are nonrepetitive on walks, and provide a conjecture that would imply that every graph with maximum degree $\Delta$ has a $f(\Delta)$-colouring that is nonrepetitive on walks. We prove that every graph with treewidth $k$ and maximum degree $\Delta$ has a $O(k\Delta)$-colouring that is nonrepetitive on paths, and a $O(k\Delta^3)$-colouring that is nonrepetitive on walks.
\end{abstract}

\maketitle

\sect{Introduction}{Introduction}
%%%%%%%%%%%%%%%%%%%%%%%%%%%%%%%%%%%%%%%%%%%%%%%%%%%%%%%%%%%%%%%%%%%

We consider simple, finite, undirected graphs $G$ with vertex set $V(G)$, edge set $E(G)$, and maximum degree $\Delta(G)$. Let $[t]:=\{1,2,\dots,t\}$. A \emph{walk} in $G$ is a sequence $v_1,v_2,\dots,v_t$ of vertices of $G$, such that $v_iv_{i+1}\in E(G)$ for all $i\in[t-1]$. A $k$-\emph{colouring} of $G$ is a function $f$ that assigns one of $k$ colours to each vertex of $G$. A walk $v_1,v_2,\dots,v_{2t}$ is \emph{repetitively} coloured by $f$ if $f(v_i)=f(v_{t+i})$ for all $i\in[t]$. A walk $v_1,v_2,\dots,v_{2t}$ is \emph{boring} if $v_i=v_{t+i}$ for all $i\in[t]$. Of course, a boring walk is repetitively coloured by every colouring. We say a colouring $f$ is \emph{nonrepetitive on walks} (or \emph{walk-nonrepetitive}) if the only walks that are repetitively coloured by $f$ are boring. Let $\sigma(G)$ denote the minimum $k$ such that $G$ has a $k$-colouring that is nonrepetitive on walks. 

A walk $v_1,v_2,\dots,v_t$ is a \emph{path} if $v_i\ne v_j$ for all distinct $i,j\in[t]$. A colouring $f$ is \emph{nonrepetitive on paths} (or \emph{path-nonrepetitive}) if no path of $G$ is repetitively coloured by $f$. Let $\pi(G)$ denote the minimum $k$ such that $G$ has a $k$-colouring that is nonrepetitive on paths. Observe that a colouring that is path-nonrepetitive is \emph{proper}, in the sense that adjacent vertices receive distinct colours. Moreover, a path-nonrepetitive colouring has no $2$-coloured $P_4$ (a path on four vertices). A proper colouring with no $2$-coloured $P_4$ is called a \emph{star colouring} since each bichromatic subgraph is a star forest; see \citep{Albertson-EJC04, Grunbaum73, FRR-JGT04, Borodin-DM79, NesOdM-03}. The \emph{star chromatic number} $\chi_{\text{st}}(G)$ is the minimum number of colours in a proper colouring of $G$ with no $2$-coloured $P_4$. Thus
\begin{equation}
\chi(G)\leq\chi_{\text{st}}(G)\leq\pi(G)\leq\sigma(G).
\end{equation}

Path-nonrepetitive colourings are widely studied \citep{Grytczuk-EJC02, BreakingRhythm, CG-ENDM07, BGKNP-NonRepTree-DM07, Currie-TCS05, AG-EuroComb05, KP, AGHR-RSA02, BV-NonRepVertex, BK-AC04}; see the survey by \citet{Gryczuk-IJMMS07}. Nonrepetitive edge colourings have also been considered \citep{BV-NonRepEdge, AGHR-RSA02}.

The seminal result in this field is by \citet{Thue06}, who in 1906 proved\footnote{The nonrepetitive $3$-colouring of $P_n$ by \citet{Thue06} is obtained as follows. Given a nonrepetitive sequence over $\{1,2,3\}$, replace each $1$ by the sequence $12312$,  replace each $2$ by the sequence $131232$, and replace each $3$ by the sequence $1323132$. \citet{Thue06} proved that the new sequence is nonrepetitive. Thus arbitrarily long paths can be nonrepetitively $3$-coloured.} that the $n$-vertex path $P_n$ satisfies
\begin{equation}
\label{Path}
\pi(P_n)=%\sigma(P_n)=
\begin{cases}
n	&\text{ if }n\leq 2,\\
3	&\text{ otherwise.}
\end{cases}
\end{equation}
A result by \citet{KP} (see \lemref{PathLoop}) implies 
\begin{equation}
\sigma(P_n)\leq 4\enspace.
\end{equation}
\citet{Currie-EJC02} proved that the $n$-vertex cycle $C_n$ satisfies
\begin{equation}
\label{Cycle}
\pi(C_n)=
\begin{cases}
4	&\text{if }n\in\{5,7,9,10,14,17\},\\
3	&\text{otherwise.}
\end{cases}
\end{equation}

Let $\pi(\Delta)$ and $\sigma(\Delta)$ denote the maximum of $\pi(G)$ and $\sigma(G)$, taken over all graphs $G$ with maximum degree $\Delta(G)\leq\Delta$. Now $\pi(2)=4$ by \eqref{Path} and \eqref{Cycle}. In general, \citet{AGHR-RSA02} proved that 
\begin{equation}
\label{PiDegree}
\frac{\alpha\Delta^2}{\log\Delta}\leq\pi(\Delta)\leq\beta\Delta^2,
\end{equation}
for some constants $\alpha$ and $\beta$. The upper bound was proved using the 
Lov\'{a}sz Local Lemma, and the lower bound is attained by a random graph. 

In \secref{Bounded} we study whether $\sigma(\Delta)$ is finite, and provide a natural conjecture that would imply an affirmative answer. 

In \secref{Treewidth} we study path- and walk-nonrepetitive colourings of graphs of bounded treewidth\footnote{The \emph{treewidth} of a graph $G$ can be defined to be the minimum integer $k$ such that $G$ is a subgraph of a chordal graph with no clique on $k+2$ vertices. Treewidth is an important graph parameter, especially in structural graph theory and algorithmic graph theory; see the surveys \citep{Bodlaender-TCS98, Reed-AlgoTreeWidth03}.}. \citet{KP} and \citet{BV-NonRepVertex} independently proved that graphs of bounded treewidth have bounded $\pi$. The best bound is due to \citet{KP} who proved that $\pi(G)\leq4^k$ for every graph $G$ with treewidth at most $k$. Whether there is a polynomial bound on $\pi$ for graphs of treewidth $k$ is an open question. 
We answer this problem in the affirmative under the additional assumption of bounded degree. In particular, we prove a \Oh{k\Delta} upper bound on $\pi$, and a \Oh{k\Delta^3} upper bound on $\sigma$.

%\comment{DW: My Guess is that there is $k$-tree $G$ with $\pi(G)\geq c^k$ for some constant $c>1$.}

In \secref{Subdivision} we prove that every graph has a subdivision that admits a path-nonrepetitive $4$-colouring; the best previous bound was $5$. 

In \secref{Density} we determine the maximum density of a graph that admits a path-nonrepetitive $k$-colouring, and prove bounds on the maximum density for walk-nonrepetitive $k$-colourings. 

%Grytczuk-Berge, AG-ICGT05, 

\sect{Bounded}{Is $\sigma(\Delta)$ bounded?}
%%%%%%%%%%%%%%%%%%%%%%%%%%%%%%%%%%%%%%%%%%%%%%%%%%%%%%%%%%%%%%%%%%%

Consider the following elementary lower bound on $\sigma$, where $G^2$ is the \emph{square} graph of $G$. That is, $V(G^2)=V(G)$, and $vw\in E(G^2)$ if and only if the distance between $v$ and $w$ in $G$ is at most $2$. A proper colouring of $G^2$ is called a \emph{distance-$2$} colouring of $G$.

\lem{WalkLowerBound}{Every walk-nonrepetitive colouring of a graph $G$ is distance-$2$. Thus $\sigma(G)\geq\chi(G^2)\geq\Delta(G)+1$.}

\prf{Consider a walk-nonrepetitive colouring of $G$. Adjacent vertices $v$ and $w$ receive distinct colours, as otherwise $v,w$ would be a repetitively coloured path. If $u,v,w$ is a path, and $u$ and $w$ receive the same colour, then the non-boring walk $u,v,w,v$ is repetitively coloured. Thus vertices at distance at most $2$ receive distinct colours. Hence $\sigma(G)\geq\chi(G^2)$. In a distance-$2$ colouring, each vertex and its neighbours all receive distinct colours. Thus $\chi(G^2)\geq\Delta(G)+1$.}

Hence $\Delta(G)$ is a lower bound on $\sigma(G)$. Whether high degree is the only obstruction for bounded $\sigma$ is an open problem.

%\opn{sigmapi}{Is there a function $f$ such that $\sigma(G)\leq f(\pi(G),\Delta(G))$ for every graph $G$?}

\opn{Degree}{Is there a function $f$ such that $\sigma(\Delta)\leq f(\Delta)$?}

%By \eqnref{PiDegree}, a positive answer to \openref{sigmapi} would imply  a positive answer to \openref{Degree}.

First we answer \openref{Degree} in the affirmative for $\Delta=2$.
The following lemma will be useful.

\lem{useful}{Fix a distance-$2$ colouring of a graph $G$. If  $W=(v_1,v_2,\dots,v_{2t})$ is a repetitively coloured non-boring walk in $G$, then $v_i\ne v_{t+i}$ for all $i\in[t]$.}

\prf{Suppose on the contrary that $v_i=v_{t+i}$ for some $i\in[t-1]$. Since $W$ is repetitively coloured, $c(v_{i+1})=c(v_{t+i+1})$. Each neighbour of $v_i$ receives a distinct colour. Thus $v_{i+1}=v_{t+i+1}$. By induction, $v_j=v_{t+j}$ for all $j\in[i,t]$. By the same argument, $v_j=v_{t+j}$ for all $j\in[1,i]$. Thus 
$W$ is boring, which is the desired contradiction.}

\prop{PathCycle}{$\sigma(2)\leq 5$.}

\prf{A result by \citet{KP} implies that $\sigma(P_n)\leq 4$ (see \lemref{PathLoop}).
Thus it suffices to prove that $\sigma(C_n)\leq 5$. 
Fix a walk-nonrepetitive $4$-colouring of the path $(v_1,v_2,\dots,v_{2n-4})$.
Thus for some $i\in[1,n-2]$, the vertices $v_i$ and $v_{n+i-2}$ receive distinct colours.
Create a cycle $C_n$ from the sub-path $v_i,v_{i+1},\dots,v_{n+i-2}$ by adding one vertex $x$ adjacent to $v_i$ and $v_{n+i-2}$. 
Colour $x$ with a fifth colour.
Observe that since $v_i$ and $v_{n+i-2}$ receive distinct colours, 
the colouring of $C_n$ is distance-$2$.
Suppose on the contrary that $C_n$ has a repetitively coloured walk $W=y_1,y_2,\dots,y_{2t}$. 
If $x$ is not in $W$, then $W$ is a repetitively coloured walk in the starting path, which is a contradiction.
Thus $x=y_i$ for some $i\in[t]$ (with loss of generality, by considering the reverse of $W$). 
Since $x$ is the only vertex receiving the fifth colour 
and $W$ is repetitive, $x=y_{t+i}$.
By \lemref{useful}, $W$ is boring.
Hence the $5$-colouring of $C_n$ is walk-nonrepetitive.}

Below we propose a conjecture that would imply a positive answer to \openref{Degree}. First consider the following lemma which is a slight generalisation of a result by \citet{BV-NonRepEdge}. A walk $v_1,v_2,\dots,v_t$ has \emph{length} $t$ and \emph{order} $|\{v_i:1\leq i\leq t\}|$. That is, the order is the number of distinct vertices in the walk.

\lem{SmallWalks}{Suppose that in some coloured graph, there is a repetitively coloured non-boring walk. Then there is a repetitively coloured non-boring walk of order $k$ and length at most $2k^2$.}

\begin{figure}[!ht]
\begin{center}
\includegraphics{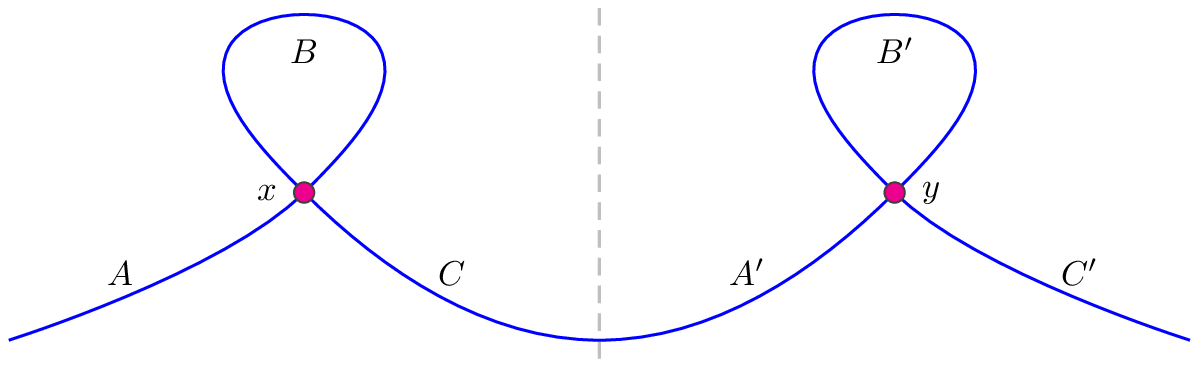}
\caption{Illustration for the proof of \lemref{SmallWalks}.}
\label{fig:Shorten} 
\end{center}
\end{figure}

\prf{Let $k$ be the minimum order of a repetitively coloured non-boring walk. 
Let $W=v_1,v_2,\dots,v_{2t}$ be a repetitively coloured non-boring walk of order $k$ and with $t$ minimum. If $2t\leq2k^2$, then we are done. Now assume that $t>k^2$. By the pigeonhole principle, there is a vertex $x$ that appears at least $k+1$ times in $v_1,v_2,\dots,v_t$. Thus there is a vertex $y$ that appears at least twice in the set $\{v_{t+i}:v_i=x,i\in[t]\}$. As illustrated in Figure~\ref{fig:Shorten}, $W=AxBxCA'yB'yC'$ for some walks $A,B,C,A',B',C'$ with $|A|=|A'|$, $|B|=|B'|$, and $|C|=|C'|$. Consider the walk $U:=AxCA'yC'$. If $U$ is not boring, then it is a repetitively coloured non-boring walk of order at most $k$ and length less than $2t$, which contradicts the minimality of $W$. Otherwise $U$ is boring, implying $x=y$, $A=A'$, and $C=C'$. Thus $B\ne B'$ since $W$ is not boring, implying $xBxB'$ is a repetitively coloured non-boring walk of order at most $k$ and length less than $2t$, which contradicts the minimality of $W$.}

We conjecture the following strengthening of \lemref{SmallWalks}.

%\con{SmallWalks}{There exists a function $h$ such that for every graph $G$ and for every path-nonrepetitive distance-$2$ colouring of $G$ with $c$ colours, if $G$ contains a repetitively coloured non-boring walk, then $G$ contains a repetitively coloured non-boring walk of order $k$ and length at most $h(c)\cdot k$.}

\con{SmallWalks}{Let $G$ be a graph. Consider a path-nonrepetitive distance-$2$ colouring of $G$ with $c$ colours, such that $G$ contains a repetitively coloured non-boring walk. Then $G$ contains a repetitively coloured non-boring walk of order $k$ and length at most $h(c)\cdot k$, for some function $h$ that only depends on $c$.}

\thm{sigmaBounded}{If \conref{SmallWalks} is true, then there is a 
function $f$ for which $\sigma(\Delta)\leq f(\Delta)$. That is, every graph $G$ has a walk-nonrepetitive colouring with $f(\Delta(G))$ colours.}

\thmref{sigmaBounded} is proved using the Lov\'{a}sz Local Lemma \citep{EL75}. 

\lem[\citep{EL75}]{LLL}{Let $\A=\A_1\cup\A_2\cup\dots\cup\A_r$ be a partition of a set of `bad' events \A. Suppose that there are sets of real numbers $\{p_i\in[0,1):i\in[r]\}$, $\{x_i\in[0,1):i\in[r]\}$, and $\{D_{ij}\geq0:i,j\in[r]\}$ such that the following conditions are satisfied by every event $A\in\A_i$:
\begin{itemize}
\item the probability $\displaystyle
\Prob(A)\;\leq\;p_i\;\leq\;x_i\cdot\!\!\!\!\prod_{j=1}^r(1-x_j)^{D_{ij}}\enspace$, and
\item $A$ is mutually independent of $\A\setminus(\{A\}\cup\D_A)$, for some $\D_A\subseteq\A$ with $|\D_A\cap\A_j|\leq D_{ij}$ for all $j\in[r]$. 
\end{itemize}
Then 
\begin{equation*}
\Prob\bracket{\bigwedge_{A\in\A} \overline{A}}\;\geq\;\prod_{i=1}^r(1-x_i)^{|\A_i|}\;>\;0\enspace.
\end{equation*}
That is, with positive probability, no event in $\A$ occurs.}

\prf[Proof of \thmref{sigmaBounded}]{Let $f_1$ be a path-nonrepetitive colouring of $G$ with $\pi(G)$ colours. Let $f_2$ be a distance-$2$ colouring of $G$ with $\chi(G^2)$ colours. Note that $\pi(G)\leq \beta\Delta^2$ for some constant $\beta$ by Equation~\eqref{PiDegree}, and $\chi(G^2)\leq\Delta(G^2)+1\leq\Delta^2+1$ by a greedy colouring of $G^2$. Hence $f_1$ and $f_2$ together define a path-nonrepetitive distance-$2$ colouring of $G$. The number of colours $\pi(G)\cdot\chi(G^2)$ is bounded by a function solely of $\Delta(G)$. Consider this initial colouring to be fixed. Let $c$ be a positive integer to be specified later. For each vertex $v$ of $G$, choose a third colour $f_3(v)\in[c]$ independently and randomly. Let $f$ be the colouring defined by $f(v)=(f_1(v),f_2(v),f_3(v))$ for all vertices $v$.

Let $h:=h(\pi(G)\cdot\chi(G^2))$ from \conref{SmallWalks}. A non-boring walk $v_1,v_2,\dots,v_{2t}$ of order $i$ is \emph{interesting} if its length $2t\leq hi$, and $f_1(v_j)=f_1(v_{t+j})$ and $f_2(v_j)=f_2(v_{t+j})$ for all $j\in[t]$. For each interesting walk $W$, let $A_W$ be the event that $W$ is repetitively coloured by $f$. Let $\A_i$ be the set of events $A_W$, where $W$ is an interesting walk of order $i$. Let $\A=\bigcup_i\A_i$. 

We will apply \lemref{LLL} to prove that, with positive probability, no event $A_W$ occurs. This will imply that there exists a colouring $f_3$ such that no interesting walk is repetitively coloured by $f$. A non-boring non-interesting walk $v_1,v_2,\dots,v_{2t}$ of order $i$ satisfies (a) $2t>hi$, or (b) $f_1(v_j)\ne f_1(v_{t+j})$ or $f_2(v_j)\ne f_2(v_{t+j})$ for some $j\in[t]$. In case (a), by the assumed truth of \conref{SmallWalks}, $W$ is not repetitively coloured by $f$. In case (b), $f(v_j)\ne f(v_{t+j})$ and $W$ is not repetitively coloured by $f$. Thus no non-boring walk is repetitively coloured by $f$, as desired.

Consider an interesting walk $W=v_1,v_2,\dots,v_{2t}$ of order $i$. 

We claim that $v_\ell\ne v_{t+\ell}$ for all $\ell\in[t]$. 
Suppose on the contrary that $v_\ell=v_{t+\ell}$ for some $\ell\in[t]$. 
Since $W$ is not boring, $v_j\ne v_{t+j}$ for some $j\in[t]$. 
Thus $v_j=v_{t+j}$ and $v_{j+1}\ne v_{t+j+1}$  for some $j\in[t]$ (where $v_{t+t+1}$ means $v_1$). Since $W$ is interesting, $f_2(v_{j+1})=f_2(v_{t+j+1})$, which is a contradiction since $v_{j+1}$ and $v_{t+j+1}$ have a common neighbour $v_j$ ($=v_{t+j}$). Thus $v_j\ne v_{t+j}$ for all $j\in[t]$, as claimed. 

This claim implies that for each of the $i$ vertices $x$ in $W$, there is at least one other vertex $y$ in $W$, such that $f_3(x)=f_3(y)$ must hold for $W$ to be repetitively coloured. Hence at most $c^{i/2}$ of the $c^i$ possible colourings of $W$ under $f_3$, lead to repetitive colourings of $W$ under $f$. Thus the probability $\Prob(A_W)\leq p_i:=c^{-i/2}$, and 
\lemref{LLL} can be applied as long as
\begin{equation}
\label{eqn:One}
c^{-i/2}\;\leq\;x_i\cdot\!\!\!\!\prod_j(1-x_j)^{D_{ij}}\enspace,
\end{equation}

Every vertex is in at most $hj\Delta^{hj}$ interesting walks of order $j$. Thus an interesting walk of order $i$ shares a vertex with at most $hij\Delta^{hj}$ interesting walks of order $j$. Thus we can take $D_{ij}:=hij\Delta^{hj}$. 
Define $x_i:=(2\Delta^{h})^{-i}$. Note that $x_i\leq\half$. So $1-x_i\geq\e^{-2x_i}$. Thus to prove \eqref{eqn:One} it suffices to prove that 
\begin{align*}
&
c^{-i/2}
\;\leq\;
x_i\cdot\!\!\!\!\prod_j\e^{-2x_jD_{ij}}\enspace,\\
\Longleftarrow\;\;\;&
c^{-i/2}
\;\leq\;
(2\Delta^{h})^{-i}\cdot\!\!\!\!\prod_j\e^{-2(2\Delta^{h})^{-j}hij\Delta^{hj}}\enspace,\\
\Longleftarrow\;\;\;&
c^{-1/2}
\;\leq\;
(2\Delta^{h})^{-1}\cdot\!\!\!\!\prod_j\e^{-2(2)^{-j}hj}\enspace,\\
\Longleftarrow\;\;\;&
c^{-1/2}
\;\leq\;
(2\Delta^{h})^{-1}\e^{-2h\sum_jj2^{-j}}\enspace,\\
\Longleftarrow\;\;\;&
c^{-1/2}
\;\leq\;
(2\Delta^{h})^{-1}\e^{-4h}\enspace,\\
\Longleftarrow\;\;\;&
%c^{1/2}
%&\;\geq\;
%2(\e^4\Delta)^h\enspace,\\
c\;\;\;\;\;\;\;\;\geq\;
4(\e^4\Delta)^{2h}\enspace.
\end{align*}
Choose $c$ to be the minimum integer that satisfies this inequality, and the lemma is applicable. We obtain a $c$-colouring $f_3$ of $G$ such that $f$ is nonrepetitive on walks. The number of colours in $f$ is at most $h\ceil{4(\e^4\Delta)^{2h}}$, which is a function solely of $\Delta$.}

\sect{Treewidth}{Trees and Treewidth}
%%%%%%%%%%%%%%%%%%%%%%%%%%%%%%%%%%%%%%%%%%%%%%%%%%%%%%%%%%%%%%%%%%%

We start this section by considering walk-nonrepetitive colourings of trees.

\thm{Trees}{Let $T$ be a tree. A colouring $c$ of $T$ is walk-nonrepetitive if and only if $c$ is path-nonrepetitive and distance-$2$.}

\prf{For every graph, every walk-nonrepetitive colouring is path-nonrepetitive (by definition) and distance-$2$ (by \lemref{WalkLowerBound}). 

Now fix a path-nonrepetitive distance-$2$ colouring $c$ of $T$. Suppose on the contrary that $T$ has a repetitively coloured non-boring walk. Let $W=(v_1,v_2,\dots,v_{2t})$ be a repetitively coloured non-boring walk in $T$ of minimum length. Some vertex is repeated in $W$, as otherwise $W$ would be a repetitively coloured path. By considering the reverse of $W$, without loss of generality, $v_i=v_j$ for some $i\in[1,t-1]$ and $j\in[i+2,2t]$. Choose $i$ and $j$ to minimise $j-i$. Thus $v_i$ is not in the sub-walk $(v_{i+1},v_{i+2},\dots,v_{j-1})$. Since $T$ is a tree, $v_{i+1}=v_{j-1}$. Thus $i+1=j-1$, as otherwise $j-i$ is not minimised. That is, $v_i=v_{i+2}$. Assuming $i\neq t-1$, since $W$ is repetitively coloured, $c(v_{t+i})=c(v_{t+i+2})$, which implies that $v_{t+i}=v_{t+i+2}$ because $c$ is a distance-$2$ colouring. Thus, even if $i=t-1$, deleting the vertices $v_{i},v_{i+1},v_{t+i},v_{t+i+1}$ from $W$, gives a walk $(v_1,v_2,\dots,v_{i-1},v_{i+2},\dots,v_{t+i-1},v_{t+i+2},\dots,v_{2t})$ that is also repetitively coloured. This contradicts the minimality of the length of $W$.}

Note that \thmref{Trees} implies that \conref{SmallWalks} is vacuously true for trees.

Since every tree $T$ has a path-nonrepetitive $4$-colouring \citep{KP} and a distance-$2$ $(\Delta(T)+1)$-colouring, \thmref{Trees} implies the following result, where the lower bound is \lemref{WalkLowerBound}.

\cor{WalkTree}{Every tree $T$ satisfies $\Delta(T)+1\leq\sigma(T)\leq 4(\Delta(T)+1)$.}

%The \emph{treewidth} of a graph $G$ can be defined to be the minimum integer $k$ such that $G$ is a subgraph of a chordal graph with no clique on $k+2$ vertices. Treewidth is an important graph parameter, especially in structural graph theory and algorithmic graph theory; see the surveys \citep{Bodlaender-TCS98, Reed-AlgoTreeWidth03}. 

%\citet{KP} and \citet{BV-NonRepVertex} independently proved that graphs of bounded treewidth have bounded $\pi$. The best bound is due to \citet{KP} who proved that $\pi(G)\leq4^k$ for every graph $G$ with treewidth at most $k$. Whether there is a polynomial bound on $\pi$ for graphs of treewidth $k$ is an open question. 

In the remainder of this section we prove the following polynomial upper bounds on $\pi$ and $\sigma$ in terms of the treewidth and maximum degree of a graph.

\thm{TreewidthDegree}{Every graph $G$ with treewidth $k$ and maximum degree $\Delta\geq1$ satisfies $\pi(G)\leq ck\Delta$ and $\sigma(G)\leq ck\Delta^3$ for some constant $c$.}

We prove \thmref{TreewidthDegree} by a series of lemmas. The first is by \citet{KP}\footnote{The $4$-colouring in \lemref{PathLoop} is obtained as follows. Given a nonrepetitive sequence on $\{1,2,3\}$, insert the symbol $4$ between consecutive block of length two. For example, from the sequence $123132123$ we obtain $1243143241243$.}. 

\lem[\citep{KP}]{PathLoop}{Let $P^+$ be the pseudograph obtained from a path $P$ by adding a loop at each vertex. Then $\sigma(P^+)\leq 4$.}

%\comment{DW: The proof of \lemref{PathLoop} implies that if a path $P$ is coloured so that it is nonrepetitive on paths and is a proper colouring of $P^2$, then it is nonrepetitive on walks. For what other graphs can we generalise this statement, possibly with other conditions on the colouring?
%JB: I think the idea is that a proper  coloring of $G^2$ kills the walks of length 3 on 3  vertices. I will think about it.}

Now we introduce some definitions by \citet{KP}. A \emph{levelling} of a graph $G$ is a function $\lambda:V(G)\rightarrow\Z$ such that $|\lambda(v)-\lambda(w)|\leq1$ for every edge $vw\in E(G)$. Let $G_{\lambda=k}$ and $G_{\lambda>k}$ denote the subgraphs of $G$ respectively induced by $\{v\in V(G):\lambda(v)=k\}$ and $\{v\in V(G):\lambda(v)>k\}$. The \emph{$k$-shadow} of a subgraph $H$ of $G$ is the set of vertices in $G_{\lambda=k}$ adjacent to some vertex in $H$. A levelling $\lambda$ is \emph{shadow-complete} if the $k$-shadow of every component of $G_{\lambda>k}$ induces a clique. \citet{KP} proved the following lemma for repetitively coloured paths. We show that the same proof works for repetitively coloured walks. 

\lem{LevellingWalks}{For every levelling $\lambda$ of a graph $G$, there is a $4$-colouring of $G$, such that every repetitively coloured walk $v_1,v_2,\dots,v_{2t}$ satisfies $\lambda(v_j)=\lambda(v_{t+j})$ for all $j\in[t]$.}

\prf{The levelling $\lambda$ can be thought of as a homomorphism from $G$ into $P^+$, for some path $P$. By \lemref{PathLoop}, $P^+$ has a 4-colouring that is nonrepetitive on walks. Colour each vertex $v$ of $G$ by the colour assigned to $\lambda(v)$ (thought of as a vertex of $P^+$). Suppose $v_1,v_2,\dots,v_{2t}$ is a repetitively coloured walk in $G$. Thus $\lambda(v_1),\lambda(v_2),\dots,\lambda(v_{2t})$ is a repetitively coloured walk in $P^+$. Since the 4-colouring of $P^+$ is nonrepetitive on walks, $\lambda(v_1),\lambda(v_2),\dots,\lambda(v_{2t})$ is boring. That is, $\lambda(v_j)=\lambda(v_{t+j})$ for all $j\in[t]$.}

\lem[\citep{KP}]{Shadow}{If $\lambda$ is a shadow-complete levelling of a graph $G$, then $$\pi(G)\leq 4\cdot\max_k\pi(G_{\lambda=k}).$$}

Now we generalise \lemref{Shadow} for walks. 

\lem{ShadowWalks}{If $H$ is a subgraph of a graph $G$, and $\lambda$ is a shadow-complete levelling of $G$, then 
$$\sigma(H)\;\leq\; 4\,\chi(H^2)\cdot\max_k\sigma(G_{\lambda=k})
\;\leq\; 4(\Delta(H)^2+1)\cdot\max_k\sigma(G_{\lambda=k}).$$}

\prf{Let $c_1$ be the $4$-colouring of $G$ from \lemref{LevellingWalks}. Thus every repetitively coloured walk $v_1,v_2,\dots,v_{2t}$ satisfies $\lambda(v_j)=\lambda(v_{t+j})$ for all $j\in[t]$. Let $c_2$ be an optimal walk-nonrepetitive colouring of each level $G_{\lambda=k}$. Let $c_3$ be a proper $\chi(H^2)$-colouring of $H^2$. The second inequality in the lemma follows from the first since $\chi(H^2)\leq\Delta(H)^2+1$. Let $c(v):=(c_1(v),c_2(v),c_3(v))$ for each vertex $v$ of $H$. We claim that $c$ is nonrepetitive on walks in $H$.

Suppose on the contrary that $W=v_1,\dots,v_{2t}$ is a non-boring walk in $H$ that is repetitively coloured by $c$. Then $W$ is repetitively coloured by each of $c_1$, $c_2$, and $c_3$. Thus $\lambda(v_i)=\lambda(v_{t+i})$ for all $i\in[t]$ by \lemref{LevellingWalks}. Let $W_k$ be the sequence (allowing  repetitions) of vertices $v_i\in W$ such that $\lambda(v_i)=k$. Since $v_i\in W_k$ if and only if $v_{t+i}\in W_k$, each sequence $W_k$ is repetitively coloured. That is, if $W_k=x_1,\dots,x_{2s}$ then $c(x_i)=c(x_{s+i})$ for all $i\in[s]$. 

Let $k$ be the minimum level containing a vertex in $W$. Let $v_i$ and $v_j$ be consecutive vertices in $W_k$ with $i<j$. If $j=i+1$ then $v_iv_j$ is an edge of $W$. Otherwise there is walk from $v_i$ to $v_j$ in $G_{\lambda>k}$ (since $k$ was chosen minimum), implying $v_iv_j$ is an edge of $G$ (since $\lambda$ is shadow-complete). Thus $W_k$ forms a walk in $G_{\lambda=k}$ that is repetitively coloured by $c_2$. Hence $W_k$ is boring. In particular, some vertex $v_i=v_{t+i}$ is in $W_k$. Since $W$ is not boring, $v_j\ne v_{t+j}$ for some $j\in[t]$. Without loss of generality, $i<j$ and $v_\ell=v_{t+\ell}$ for all $\ell\in[i,j-1]$. Thus $v_j$ and $v_{t+j}$ have a common neighbour $v_{j-1}=v_{t+j-1}$ in $H$, which implies that $c_3(v_j)\ne c_3(v_{t+j})$. But $c(v_j)=c(v_{t+j})$ since $W$ is repetitively coloured, which is the desired contradiction.}

Note that some dependence on $\Delta(H)$ in \lemref{ShadowWalks} is unavoidable, since $\sigma(H)\geq\chi(H^2)\geq\Delta(H)+1$. 

\lemref{ShadowWalks} enables the following strengtening of \corref{WalkTree}.

\lem{TreeWalks}{Every tree $T$ satisfies $\Delta(T)+1\leq\sigma(T)\leq 4\,\Delta(T)$.}

\prf{Let $r$ be a leaf vertex of $T$. Let $\lambda(v)$ be the distance from $r$ to $v$ in $T$. Then $\lambda$ is a shadow-complete levelling of $T$ in which each level is an independent set. A greedy algorithm proves that $\chi(T^2)\leq\Delta(T)+1$. Thus \lemref{ShadowWalks} implies that $\sigma(T)\leq 4\,\Delta(T)+4$. Observe that the proof of \lemref{ShadowWalks} only requires $c_3(v)\ne c_3(w)$ whenever $v$ and $w$ are in the same level and have a common parent. Since $r$ is a leaf, each vertex has at most $\Delta(T)-1$ children. Thus a greedy algorithm produces a $\Delta(T)$-colouring with this property. Hence $\sigma(T)\leq4\,\Delta(T)$.}

A \emph{tree-partition} of a graph $G$ is a partition of its vertices into sets (called \emph{bags}) such that the graph obtained from $G$ by identifying the vertices in each bag is a forest (after deleting loops and replacing parallel edges by a single edge)\footnote{The proof by \citet{KP} that $\pi(G)\leq4^k$ for graphs with treewidth at most $k$ can also be described using tree-partitions; cf.~\citep{DMW-SJC05}.}.
 
\lem{TreePartitionLevelling}{Let $G$ be a graph with a tree-partition in which every bag has at most $\ell$ vertices. Then $G$ is a subgraph of a graph $G'$ that has a shadow-complete levelling in which each level satisfies $$\pi(G'_{\lambda=k})\leq\sigma(G'_{\lambda=k})\leq\ell.$$}

\prf{Let $G'$ be the graph obtained from $G$ by adding an edge between all pairs of nonadjacent vertices in a common bag. Let $F$ be the forest obtained from $G'$ by identifying the vertices in each bag. Root each component of $F$. Consider a vertex $v$ of $G'$ that is in the bag that corresponds to node $x$ of $F$. Let $\lambda(v)$ be the distance between $x$ and the root of the tree component of $F$ that contains $x$. Clearly $\lambda$ is a levelling of $G'$. The $k$-shadow of each connected component of $G'_{\lambda>k}$ is contained in a single bag, and thus induces a clique on at most $\ell$ vertices. Hence $\lambda$ is shadow-complete. By colouring the vertices within each bag with distinct colors, we have $\pi(G'_{\lambda=k})\leq\sigma(G'_{\lambda=k})\leq\ell$.}

\threelemref{Shadow}{ShadowWalks}{TreePartitionLevelling} imply:

\lem{TreePartitionColouring}{If a graph $G$ has a tree-partition in which every bag has at most $\ell$ vertices, then $\pi(G)\leq 4\ell$ and $\sigma(G)\leq4\ell(\Delta(G)^2+1)$.}

\citet{Wood-TreePartitions} proved\footnote{The proof by \citet{Wood-TreePartitions} is a minor improvement to a similar result by an anonymous referee of the paper by \citet{DO-JGT95}.} that every graph with treewidth $k$ and maximum degree $\Delta\geq1$ has a tree-partition in which every bag has at most $\frac{5}{2}(k+1)(\frac{7}{2}\Delta-1)$ vertices. With \lemref{TreePartitionColouring} this proves the following quantitative version of \thmref{TreewidthDegree}.

\thm{TreewidthDegreeQuantative}{Every graph $G$ with treewidth $k$ and maximum degree $\Delta\geq1$ satisfies $\pi(G)\leq 10(k+1)(\frac{7}{2}\Delta-1)$ and $\sigma(G)\leq 10(k+1)(\frac{7}{2}\Delta-1)(\Delta^2+1)$.}

\sect{Subdivision}{Subdivisions}
%%%%%%%%%%%%%%%%%%%%%%%%%%%%%%%%%%%%%%%%%%%%%%%%%%%%%%%%%%%%%%%%%%%

The results of \citet{Thue06} and \citet{Currie-EJC02} imply that every path and every cycle has a subdivision $H$ with $\pi(H)=3$. \citet{BGKNP-NonRepTree-DM07} proved that every tree has a subdivision $H$ such that $\pi(H)=3$. Which graphs have a subdivision $H$ with $\pi(H)=3$ is an open problem \citep{Gryczuk-IJMMS07}. \citet{Gryczuk-IJMMS07} proved that every graph has a subdivision $H$ with $\pi(H)\leq 5$. Here we improve this bound as follows.

\thm{Subdivision}{Every graph $G$ has a subdivision $H$ with $\pi(H)\leq 4$.}

%The proof of \thmref{Subdivision} depends on the following lemma of \citet{KP}. 

%\lem[\citep{KP}]{Levelling}{For every levelling $\lambda$ of a graph $G$, there is $4$-colouring of $G$, such that every repetitively coloured path $v_1,v_2,\dots,v_{2t}$ satisfies $\lambda(v_j)=\lambda(v_{t+j})$ for all $j\in[t]$.}

\begin{figure}[!ht]
\begin{center}
\includegraphics{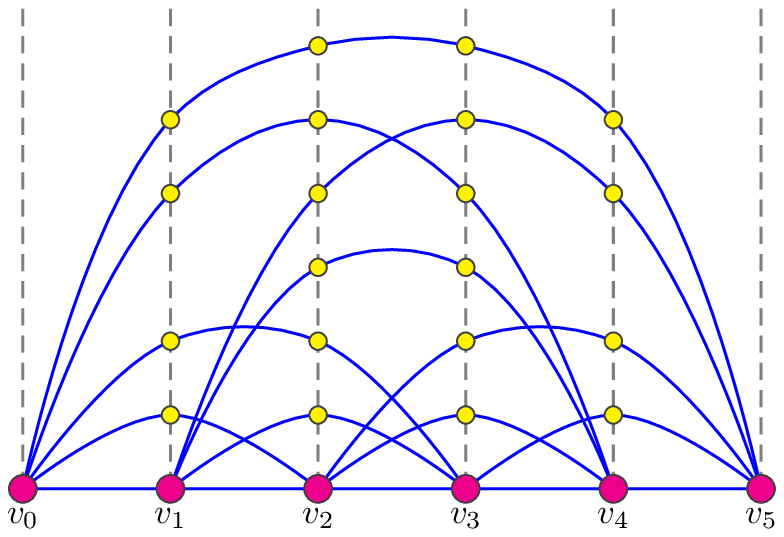}
\caption{\label{fig:Subdiv}The subdivision $H$ with $G=K_6$.}
\end{center}
\end{figure}

\begin{proof} %[Proof of \thmref{Subdivision}]
Without loss of generality $G$ is connected. Say $V(G)=\{v_0,v_1,\dots,v_{n-1}\}$. As illustrated in Figure~\ref{fig:Subdiv}, let $H$ be the subdivision of $G$ obtained by subdividing every edge $v_iv_j\in E(G)$ (with $i<j$) $j-i-1$ times. The distance of every vertex in $H$ from $v_0$ defines a levelling of $H$ such that the endpoints of every edge are in consecutive levels. By \lemref{LevellingWalks}, there is a $4$-colouring of $H$, such that for every repetitively coloured path $x_1,x_2,\dots,x_t,y_1,y_2,\dots,y_t$ in $H$, $x_j$ and $y_j$ have the same level for all $j\in[t]$. Hence there is some $j$ such that $x_{j-1}$ and $x_{j+1}$ are at the same level. Thus $x_j$ is an original vertex $v_i$ of $G$. Without loss of generality $x_{j-1}$ and $x_{j+1}$ are at level $i-1$. There is only one original vertex at level $i$. Thus $y_j$, which is also at level $i$, is a division vertex. Now $y_j$ has two neighbours in $H$, which are at levels $i-1$ and $i+1$. Thus $y_{j-1}$ and $y_{j+1}$ are at levels $i-1$ and $i+1$, which contradicts the fact that $x_{j-1}$ and $x_{j+1}$ are both at level $i-1$. Hence we have a $4$-colouring of $H$ that is nonrepetitive on paths.
\end{proof}

It is possible that every graph has a subdivision $H$ with $\pi(H)\leq3$. If true, this would provide a striking generalisation of the result of \citet{Thue06} discussed in \secref{Introduction}. 

%\comment{DW: Conjecture: Every graph $G$ has a subdivision $H$ with $\sigma(H)\leq\Delta(G)+\Oh{1}$. 
%JB:  The idea might be that we subdivide each edge $10^{10}$ times and put a palindrome- and  square-free sequence on the new edges using 4  extra symbols.
 % DW: Would it help to subdivide each edge a distinct number of times?}

%Is this related to $\chi'(G')$?

\sect{Density}{Maximum Density}
%%%%%%%%%%%%%%%%%%%%%%%%%%%%%%%%%%%%%%%%%%%%%%%%%%%%%%%%%%%%%%%%%%%

In this section we study the maximum number of edges in a nonrepetitively coloured graph.

\prop{Extremal}{The maximum number of edges in an $n$-vertex graph $G$ with $\pi(G)\leq c$ is $(c-1)n-\binom{c}{2}$.}

\prf{Say $G$ is an $n$-vertex graph with $\pi(G)\leq c$. Fix a $c$-colouring of $G$ that is nonrepetitive on paths. Say there are $n_i$ vertices in the $i$-th colour class. Every cycle receives at least three colours. Thus the subgraph induced by the vertices coloured $i$ and $j$ is a forest, and has at most $n_i+n_j-1$ edges. Hence the number of edges in $G$ is at most 
\begin{equation*}
\sum_{1\leq i<j\leq c}(n_i+n_j-1)=\sum_{1\leq i\leq c}(c-1)n_i-\binom{c}{2}
=(c-1)n-\binom{c}{2}.
\end{equation*}
This bound is attained by the graph consisting of a complete graph $K_{c-1}$ completely connected to an independent set of $n-(c-1)$ vertices, which obviously has a $c$-colouring that is nonrepetitive on paths.}

Now consider the maximum number of edges in a coloured graph that is nonrepetitive on walks. First note that the example in the proof of \propref{Extremal} is repetitive on walks. Since $\sigma(G)\geq\Delta(G)+1$ and $|E(G)|\leq\half\Delta(G)|V(G)|$, we have the trivial upper bound, 
$$|E(G)|\leq\half(\sigma(G)-1)|V(G)|.$$
This bound is tight for $\sigma=2$ (matchings) and $\sigma=3$ (cycles), but is not known to be tight for $\sigma\geq 4$. 

We have the following lower bound. 

\prop{Example}{For all $p\geq1$, there are infinitely many graphs $G$ with $\sigma(G)\leq 4p$ and $$|E(G)|\geq\eighth(3\sigma(G)-4)|V(G)|-\ninth\sigma(G)^2.$$}

\prf{Let $G$ be the lexicographic product of a path and $K_p$; that is, $G$ is the graph with a levelling $\lambda$ in which each level induces $K_p$, and every edge is present between consecutive levels. Let $c_1$ be the $4$-colouring of $G$ from \lemref{LevellingWalks}. If $v$ is the $j$-th vertex in its level, where $j\in[p]$, then let $c(v):=(c_1(v),j)$. The number of colours is $4p$. We claim that $c$ is nonrepetitive on walks in $G$. Suppose on the contrary that $W=v_1,\dots,v_{2t}$ is a non-boring walk in $G$ that is repetitively coloured by $c$. Then $W$ is repetitively coloured by $c_1$. Thus $\lambda(v_i)=\lambda(v_{t+i})$ for all $i\in[t]$ by \lemref{LevellingWalks}.
Since $W$ is not boring, some $v_i\ne v_{t+i}$. By construction, $c(v_i)\ne c(v_{t+i})$, which contradicts the assumption that $W$ is repetitively coloured. Hence $\sigma(G)\leq 4p$. Now we count the edges: 
$|E(G)|=\half(3p-1)|V(G)|-p^2$.
As a lower bound, $\sigma(G)\geq\Delta(G)+1=3p$. 
Thus $|E(G)|\geq\half(3\sigma(G)/4-1)|V(G)|-(\sigma(G)/3)^2$.}

%\opn{}{Does every graph $G$ satisfy $|E(G)|\leq(\frac{3}{8}+o(1))\sigma(G)|V(G)|$?}

\section*{Acknowledgement}
%%%%%%%%%%%%%%%%%%%%%%%%%%%%%%%%%%%%%%%%%%%%%%%%%%%%%%%%%%%%%%%%%%%%%%%%%%%%%%%%

Thanks to Carsten Thomassen whose generous hospitality at the Technical University of Denmark enabled this collaboration.

%\section*{Conjectures not for publication}
%%%%%%%%%%%%%%%%%%%%%%%%%%%%%%%%%%%%%%%%%%%%%%%%%%%%%%%%%%%%%%%%%%%%%%%%%%%%%%%%

%\comment{DW: Conjecture: The only $\Delta$-regular graph $G$ with $\sigma(G)=\Delta+1$ is $G=K_{\Delta+1}$ for all $\Delta\geq3$.

%JB: Interesting. What happens if $\Delta=3$ and the girth is large? On the other hand, it  might be again sufficient to concentrate on  walks of length 3 on 3 vertices.}

%%%%%%%%%%%%%%%%%%%%%%%%%%%%%%%%%%%%%%%%%%%%%%%%%%%%%%%%%%%%%%%%%%%%%%%%%%%%%%%%
%\bibliographystyle{myNatbibStyle}
%\bibliography{myBibliography,myConferences}
%%%%%%%%%%%%%%%%%%%%%%%%%%%%%%%%%%%%%%%%%%%%%%%%%%%%%%%%%%%%%%%%%%%%%%%%%%%%%%%%

\def\Dbar{\leavevmode\lower.6ex\hbox to 0pt{\hskip-.23ex \accent"16\hss}D}

\end{document}